\documentclass[journal]{IEEEtran}
\usepackage{multirow}
\usepackage{multicol}
\usepackage{lipsum}
\usepackage{graphicx}
\usepackage{graphics}
\usepackage{amssymb}
\usepackage{amsmath}
\usepackage{amsbsy}
\usepackage{blindtext}
\usepackage{tcolorbox}
\usepackage{amssymb}
\usepackage{scalerel}
\usepackage{mathabx}
\usepackage{verbatim}
\usepackage{booktabs}
\usepackage{rotating,tabularx}
\usepackage{mathrsfs} 
\usepackage{caption}
\usepackage{hyperref}

\usepackage{adjustbox}
\usepackage{soul}
\usepackage{xcolor}
\usepackage{cite}
\usepackage{tikz}

\usepackage{makecell}
\usepackage{color, colortbl}
\usepackage[first=0,last=9]{lcg}

\usepackage{amsmath}
\usepackage{algorithm}
\usepackage{algpseudocode}
\usepackage{amssymb}
\usepackage{geometry}
\definecolor{LightGray}{rgb}{0.7,0.7,0.7}

\newcommand*{\rn}{\textcolor{black}}

\usetikzlibrary{shapes.multipart}

\usepackage{array}
\usepackage{makecell}

\allowdisplaybreaks

\usepackage{amsthm}
\theoremstyle{definition}

\theoremstyle{remark}

\usepackage[utf8]{inputenc}
\usepackage[english]{babel}

\usepackage[caption=false,font=footnotesize]{subfig}
\usepackage[T1]{fontenc}
\usepackage{scalerel,stackengine}
\newcommand\reallywidecheck[1]{%
\savestack{\tmpbox}{\stretchto{%
  \scaleto{%
    \scalerel*[\widthof{\ensuremath{#1}}]{\kern-.6pt\bigwedge\kern-.6pt}%
    {\rule[-\textheight/2]{1ex}{\textheight}}
  }{\textheight}%
}{0.5ex}}%
\stackon[1pt]{#1}{\scalebox{-1}{\tmpbox}}%
}
\stackMath
\IEEEoverridecommandlockouts
\IEEEoverridecommandlockouts

\usepackage{mathptmx}

\newif\ifarxiv
\arxivtrue
\arxivfalse

\usepackage{enumitem}

\renewcommand{\baselinestretch}{.97}

\begin{document}

\title{\rn{Scalable Global Optimization for AC-OPF via Quadratic Convex Relaxation and Branch-and-Bound}}

\author{ Mohammadreza Iranpour$^{\ast}$, Mohammad Rasoul Narimani$^{\ast}$
\thanks{${\ast}$: Department of Electrical and Computer Engineering, California State University Northridge (CSUN). Rasoul.narimani@csun.edu. Support from NSF contract \#2308498.}%
}

\maketitle

\begin{abstract} 
\rn{The Optimal Power Flow (OPF) problem is central to the reliable and efficient operation of power systems, yet its non-convex nature poses significant challenges for finding globally optimal solutions. While convex relaxation techniques such as Quadratic Convex (QC) relaxation have shown promise in providing tight lower bounds, they typically do not guarantee global optimality. Conversely, global optimization methods like the Branch and Bound (B\&B) algorithm can ensure optimality but often suffer from high computational costs due to the large search space involved. This paper proposes a novel B\&B-assisted QC relaxation framework for solving the AC-OPF problem that leverages the strengths of both approaches. The method systematically partitions the domains of key OPF variables, specifically, voltage magnitudes and voltage angle differences, into two equal subintervals at each iteration. The QC relaxation is then applied to each subregion to compute a valid lower bound. These bounds are compared against an upper bound obtained from a feasible AC-OPF solution identified at the outset. Subregions that yield lower bounds exceeding the upper bound are pruned from the search, eliminating non-promising portions of the feasible space. By integrating the efficiency of the QC relaxation with the global search structure of the B\&B algorithm, the proposed method significantly reduces the number of subproblems explored while preserving the potential to reach the global optimum. The algorithm is implemented using the PowerModels.jl package and evaluated on a range of PGLib-OPF benchmark cases. Results demonstrate that this hybrid strategy improves computational tractability and solution quality, particularly for large OPF instances.}

\end{abstract}

\section{Introduction}
\label{sec:Introduction}

\rn{The Optimal Power Flow (OPF) problem is fundamental to the efficient and reliable operation of electric power systems. It aims to determine the optimal setpoints for system control variables such that a specified objective function, typically minimizing generation cost or system losses, is achieved. This optimization is subject to a set of nonlinear equality and inequality constraints, including the full AC power flow equations, generator capacity limits, voltage magnitude bounds at buses, thermal loss limits, and transmission line flow limits. By ensuring that these physical and operational constraints are satisfied, OPF provides feasible and economically optimal operating conditions for the power network~\cite{molzahn2019survey}. The OPF problem was first introduced by Carpentier in the 1960s~\cite{carpentier1962contribution}. Since then, it has become an important area of study in power system engineering, using many different methods, from mathematical and data-driven techniques to both deterministic and probabilistic approaches.}

\rn{Accurately modeling power flow using the nonlinear AC power flow equations leads to the AC OPF problem, which is non-convex, can exhibit multiple local optima~\cite{bukhsh2013local}, and is generally classified as NP-hard~\cite{bienstock2019strong}. Since first being formulated by Carpentier~\cite{carpentier1962contribution}, a wide variety of optimization algorithms have been applied to the OPF problem~\cite{opf_litreview1993IandII,ferc4}. Much of this research has focused on algorithms for obtaining locally optimal or approximate OPF solutions. To address these challenges, researchers have explored various solution approaches, including relaxation, approximation, and nonlinear optimization techniques. Relaxation methods, in particular, work by embedding the original non-convex feasible region defined by the power flow equations into a larger convex set, making the problem easier to solve. Transforming the OPF problem into a convex form allows the use of well-established optimization theories and algorithms, significantly improving computational tractability. In contrast, approximation methods reduce complexity by simplifying the power flow equations through approximations. While this decreases the computational burden, it also introduces inaccuracies that may limit the model’s ability to fully capture the true behavior of the power system~\cite{molzahn2019survey}.}

\rn{Relaxation methods offer the advantage of providing bounds on the optimal objective value of the original non-convex problem. In some cases, they also offer sufficient conditions to certify either infeasibility or global optimality for specific classes of power system optimization problems. This makes them particularly valuable for verifying solution quality or determining that no feasible solution exists. However, relaxation techniques can sometimes yield inaccurate results, as the relaxed problem may fail to capture all the constraints of the original formulation. As a result, solutions that appear feasible in the relaxed space may not be feasible in the actual, non-convex problem. In contrast, approximation methods, typically based on linearizing the complex power flow equations, are computationally efficient and well, suited for real-time applications and large-scale systems. However, their primary limitation is reduced accuracy, particularly in scenarios where the nonlinear characteristics of power flow are significant. Additionally, unlike relaxation techniques, approximation methods do not provide theoretical guarantees such as bounds on the optimal solution or conditions for infeasibility. A comprehensive analysis of relaxation and approximation methods is provided in \cite{molzahn2019survey}, demonstrating how these approaches can be effectively adapted to solve OPF problems, while also outlining the respective advantages and limitations of each technique.}

\rn{Branch and Bound (B\&B) is a widely used global optimization method for solving complex problems like nonconvex nonlinear programs (NLPs), mixed-integer nonlinear programs (MINLPs)~\cite{coffrin2016strengthening, hijazi2017convex}. The key idea is to break the original problem into smaller, easier-to-handle subproblems, or branches, by splitting the domain of a selected variable, often one contributing significantly to nonconvexity or feasibility violations, into two or more disjoint intervals. Each branch defines a smaller subregion of the original problem with tighter variable bounds~\cite{liberti2012branch}.
At each node in the B\&B tree, the algorithm solves a subproblem within a limited region. A node represents one of these subproblems, each with its own constraints and variable bounds. The tree is the hierarchical structure formed by these branching subproblems, representing how the original problem has been systematically divided. If the result at a node is worse than the best-known solution or if the subproblem is infeasible, it is pruned (discarded). Otherwise, it’s split further by choosing another variable to divide. This process repeats recursively, while the algorithm maintains a list of active nodes and tracks the best solution found so far. By eliminating regions that cannot contain the global optimum, B\&B progressively reduces the search space and converges to the globally optimal solution.}
\rn{Despite its strong theoretical guarantees, B\&B faces practical challenges when applied to large-scale, high-dimensional OPF problems. As the number of decision variables increases, the number of subproblems (or nodes) can grow exponentially, causing the search tree to become extremely large~\cite{burer2012non}. In OPF, these variables include continuous ones like voltage magnitudes, angles, and power flows, as well as integer variables in problems like unit commitment or network switching. Solving each of these subproblems also adds significant computational burden, especially when the models involve complex nonlinear equations or many interacting variables. As a result, even though B\&B can theoretically find the global optimum, its high computational cost has limited its practical use, particularly for large-scale or real-time power system operations.}

\rn{Convex relaxation methods improve the efficiency of the B\&B algorithm by providing lower bounds on subproblems. These bounds help determine whether a branch should be explored. If the lower bound is worse than the best known feasible solution, or if the relaxed problem is infeasible, the branch can be pruned, avoiding unnecessary computation and speeding up the solution of nonconvex problems. Relaxation methods and local solvers work well together, relaxations assess global optimality, while local solvers find good feasible solutions. Many global optimization algorithms estimate the optimality gap by comparing the relaxation’s lower bound with a local solution’s objective~\cite{xu2021verifying}. To reach a global optimum, they use strategies to reduce this gap. Stronger relaxations provide tighter bounds and more pruning, which reduces the number of subproblems explored.}

\rn{Several relaxation techniques have been developed to help solve nonconvex OPF problems, including Semidefinite Programming (SDP)\cite{lavaei2011zero}, Second-Order Cone (SOC) relaxations\cite{kocuk2018matrix}, and Quadratic Convex (QC) relaxations~\cite{hijazi2017convex, barati2020global}. SDP relaxations are known for producing very tight bounds, but they become computationally expensive as system size increases, making them less practical for large-scale or real-time use. QC relaxation has gained attention as a more scalable alternative that balances solution quality and computational cost. It works in the polar form of power flow equations and builds convex envelopes around nonconvex terms~\cite{barrows2014correcting}. While QC may not always match the tightness of SDP, it is much faster to solve and has shown strong performance in providing useful lower bounds~\cite{barati2020global, narimani2024tightening}.Recent research has focused on combining local optimization methods with relaxation techniques to efficiently find global solutions for OPF problems. Many efforts have centered on using SOC and SDP relaxations within standard spatial B\&B frameworks~\cite{gopalakrishnan2012global}. For example, \cite{gopalakrishnan2012global} presents a spatial B\&B approach for solving OPF globally, where lower bounds are computed using either the Lagrangian dual or an SDP relaxation. More recently, \cite{barati2024global} introduces two algorithms, Feasible Successive Linear Programming (FSLP) and Feasible Branch-and-Bound (FBB), that use SDP-based formulations to systematically search for globally optimal solutions.}

\rn{While SDP-based relaxations can provide tight bounds, their high computational cost limits their scalability when used with B\&B algorithms, especially for large-scale or real-time OPF problems. In contrast, QC relaxation offers a better balance between bound tightness, computational speed, and scalability~\cite{barati2020global}, making it well-suited for hybrid strategies that require both efficiency and solution quality.
This paper introduces an integration of QC relaxation within a B\&B framework to achieve global optimality for the OPF problem. We propose a QC-assisted B\&B algorithm that iteratively splits key variables, such as voltage magnitudes and phase angle differences. 
At each step, two child subproblems are generated by tightening the bounds of selected variables, and QC relaxation is applied to assess their feasibility and objective value.
First, feasibility is checked: any child that violates system constraints is immediately pruned. For the feasible children, the objective value is compared against two baselines. If it is worse than the initial QC-OPF solution, the child is pruned. Similarly, if it exceeds the objective value of the original AC-OPF problem, it is also discarded.
Only children that pass both feasibility and objective-based checks are kept for further exploration. These become the next set of parents in the B\&B tree. This strategy allows the algorithm to efficiently navigate the search space by focusing only on promising regions, while reducing unnecessary computation. By iteratively tightening variable bounds and refining the QC relaxation, the algorithm reduces the number of subproblems, improves bound accuracy, and moves toward global or near-global solutions with significantly lower computational effort.}

\rn{The rest of the paper is organized as follows: Section~\ref{sec:opf} provides a brief overview of the OPF problem. Section~\ref{sec:formulation} presents the convexification of the power flow equations using the QC relaxation method, followed by the proposed QC-Branch and Bound algorithm. In Section~\ref{sec:results}, we evaluate the method on PGLib benchmark systems and present detailed analyses. Section~\ref{sec:Conclusion} concludes the paper.}

\section{Optimal Power Flow (OPF) problem}
\label{sec:opf}

\rn{This section briefly presents the OPF formulation, a key optimization task in power system operations, based on nonlinear power flow equations. Let $\mathcal{G}$, $\mathcal{B}$, and $\mathcal{L}$ denote the sets of generators, buses, and lines, respectively. The active and reactive power generation and demand at bus $m$ are given by $P_{m,G} + j Q_{m,G}$ and $P_{m,D} + j Q_{m,D}$. Voltage magnitude and angle at bus $m\in\mathcal{B_A}$ are denoted by $V_m$ and $\theta_m$.
Each line $(m,l)\in\mathcal{L_A}$ is modeled as a $\Pi$-circuit with mutual admittance $g_{ml} + j b_{ml}$ and shunt admittance $j b_{c,ml}$. The voltage angle difference is defined as $\theta_{ml} = \theta_m - \theta_l$. Power flow variables $P_{ml}$, $Q_{ml}$, and the maximum apparent power flow limit $\overline{S}$ represent the line constraints between buses $m$ and $l$.
For each generator $i \in \mathcal{G}$, the generation cost is modeled by a quadratic function with coefficients $c_{2,i} \ge 0$, $c_{1,i}$, and $c_{0,i}$. Upper and lower bounds are denoted by $\overline{(\cdot)}$ and $\underline{(\cdot)}$. Based on these definitions, the OPF problem is formulated as follows:}


\begin{small}
\begin{subequations}
\begin{align}
&min\quad\sum_{i\in\mathcal{G}} c_{2,i}(P_{i}^{g})^2+c_{1,i}P_{i}^{g}+c_{0,i}\label{eq:obj}\\
&\nonumber \text{subject to} \quad \left(\forall i\in\mathcal{B_A}, \forall   \left(l,m\right) \in\mathcal{L_A}\right)\\
&\underline{P_{i}^g}\le P_{i}^g\le \overline{P_{i}^g}\label{eq:activelimit}\\
&\underline{Q_{i}^g}\le Q_{i}^g\le \overline{Q_{i}^g}\label{eq:reactivelimit}\\
&\underline{V_{i}}\le V_{i}\le \overline{V_{i}}\label{eq:voltagelimit}\\
&\underline{\theta_{ij}}\le \theta_{ij}\le \overline{\theta_{ij}}\label{eq:thetalimit}\\
&\theta_{ref}=0\label{eq:thetareflimit}\\
&\!\!\!\!\!\!\! g_{sh,i}\, \tilde{V}_i^2+\sum_{\substack{(l,m)\in \mathcal{L},\\\text{s.t.} \hspace{3pt} l=i}} \!\tilde{P}_{lm}+\!\!\sum_{\substack{(l,m)\in \mathcal{L},\\\text{s.t.} 
\label{eq:active_injection}\hspace{3pt} m=i}} \!\!\tilde{P}_{ml}= P_{i,G}-P_{i,D}, \\
&\!\!\!\!\!\!\! -b_{sh,i}\, \tilde{V}_i^2+\!\!\!\!\!\!\sum_{\substack{(l,m)\in \mathcal{L},\\ \text{s.t.} \hspace{3pt} l=i}} \!\!\tilde{Q}_{lm}+\!\!\!\sum_{\substack{(l,m)\in \mathcal{L},\\ \text{s.t.}\label{eq:reactive_injection}\hspace{3pt} m=i}} \!\!\!\!\tilde{Q}_{ml}=Q_{i,G}-Q_{i,D},\\
&\nonumber\!\!\!\!\!\!\! \tilde{P}_{lm} \!=\! g_{lm} \tilde{V}_l^2\! -\! g_{lm} \tilde{V}_l V_m\cos\left(\tilde{\theta}_{l}-\theta_{m}\right)\!\\
&\label{eq:qik1} \qquad\qquad -\! b_{lm} \tilde{V}_l V_m\sin\left(\tilde{\theta}_{l}-\theta_{m}\right),\\
&\!\!\!\!\!\!\! \nonumber \tilde{Q}_{lm} = -\left(b_{lm}+b_{c,lm}/2\right) \tilde{V}_l^2 + b_{lm} \tilde{V}_l V_m\cos\left(\tilde{\theta}_{l}-\theta_{m}\right)\\ 
&\!\!\!\!\!\!\! \nonumber - g_{lm} \tilde{V}_l V_m\sin\left(\tilde{\theta}_{l}-\theta_{m}\right),\\
&  \qquad\qquad (P_{lm})^{2}+(Q_{lm})^{2}\le (\overline{S_{lm}})^{2} \label{eq:thermal1}.
\end{align}
\end{subequations}
\end{small}

\rn{The AC power flow equations introduce non-convexities into the OPF problem, making it challenging to solve. These non-convexities can lead to multiple local optima, rendering the problem NP-Hard~\cite{bienstock2019strong}. To address these challenges, the QC relaxation~\cite{coffrin2015qc}, formulated in the next section, is utilized.}

\section{Problem Formulation}
\label{sec:formulation}
\rn{In this section, we present the formulation of the proposed approach. The QC relaxation method and the  B\&B algorithm are introduced and briefly explained, as they form the foundation for finding the global optimum of the OPF problem.}

\subsection{Quadratic Convex QC Relaxation formulations}

\rn{The QC relaxation~\cite{coffrin2016strengthening} is a promising approach that uses convex envelopes to approximate non-convex terms such as trigonometric functions, squared terms, and bilinear products. This technique transforms the original non-convex problem into a convex one, which is generally easier to solve. The relaxation is constructed by introducing new linear surrogate variables, $w_{ii}$, $w_{lm}$, $c_{lm}$, and $s_{lm}$, to replace nonlinear expressions. These represent, respectively, the squared voltage magnitudes, the product of voltage magnitudes at different buses, and bilinear and trilinear terms involving voltage magnitudes and trigonometric functions~\cite{coffrin2015qc}. By substituting these variables into Equation~\eqref{eq:cs}, the original non-convex problem is reformulated as a convex one.}



\vspace{-.5cm}
\begin{subequations}
\label{eq:cs}
\begin{align}
w_{ii} &= V_i^2,  & \forall i \in\mathcal{N}, \\
w_{lm} &= V_l V_m,  & \forall \left(l,m\right) \in\mathcal{L}, \\
c_{lm} & =  w_{lm} \cos\left(\theta_{lm} \right), & \forall \left(l,m\right) \in\mathcal{L}, \\
s_{lm} & = w_{lm}\sin\left(\theta_{lm} \right),  & \forall \left(l,m\right) \in\mathcal{L}.
\end{align}
\end{subequations}

\rn{For each line $\left(l,m\right)\in\mathcal{L}$, the following relationships hold among the variables $w_{ll}$, $c_{lm}$, and $s_{lm}$:}
\begin{subequations}
\label{eq:cs_relationships}
\begin{align}
\label{eq:Jabr}
&c_{lm}^2+s_{lm}^2=w_{ll}w_{mm},\\
\label{eq:cs_relationships_c}
&c_{lm}=c_{ml}, \\
\label{eq:cs_relationships_s}
&s_{lm}=-s_{ml}
\end{align}
\end{subequations}

\rn{The QC relaxation constructs convex envelopes to enclose the aforementioned non-convex terms, treating them as set-valued functions defined as follows:}

\begin{subequations}
\label{eq:product_envelopes}
\begin{align}
\label{eq:squareenvelopes}
\langle x^2\rangle^T =
\begin{cases}
\widecheck{x}: \begin{cases}\check{x} \geq x^2,\\
\widecheck{x} \leq \left({\overline{x}+\underline{x}}\right) x-{\overline{x} \underline{x}}.\\
\end{cases}
\end{cases}\\
\label{eq:mccormick}
\langle {xy}\rangle^M  =
\begin{cases}
\widecheck{xy}:\begin{cases}
\widecheck{xy} \geq {\underline{x}} y+ {\underline{y}} x-{\underline{x} \underline{y}},\\
\widecheck{xy} \geq {\overline{x}} y+ {\overline{y}} x-{\overline{x} \overline{y}},\\
\widecheck{xy} \leq {\underline{x}} y+ {\overline{y}} x-{\underline{x}} {\overline{y}},\\
\widecheck{xy} \leq {\overline{x}} y+ {\underline{y}} x-{\overline{x} \underline{y}}.\\
\end{cases}
\end{cases}
\end{align}
\end{subequations}

\rn{where $\widecheck{x}$ and $\widecheck{xy}$ are ``auxiliary'' variables representing the corresponding sets.
The envelope $\langle x^2\rangle^T$ denotes the convex hull of the squared function, while the McCormick envelope $\langle xy \rangle^M$ represents the convex hull of a bilinear product, as described in~\cite{mccormick1976}. The QC relaxation also defines convex envelopes for trigonometric functions, specifically $\left\langle \sin(x) \right\rangle^S$ and $\left\langle \cos(x) \right\rangle^C$, which provide convex approximations of the sine and cosine functions. These relaxations simplify the optimization by enclosing the non-convex trigonometric terms. The construction of these envelopes is as follows:}

\begin{subequations}
\label{eq:convex_envelopes_sin&cos}
\begin{align}
\label{eq:sine envelope}
\nonumber &\left\langle \sin(x)\right\rangle^S =\\
&\quad\;\begin{cases}
\widecheck{S}:\begin{cases}
\widecheck{S}\leq\cos\left(\frac{x^m}{2}\right)\left(x-\frac{x^m}{2}\right)+\sin \left(\frac{x^m}{2}\right),\\
\widecheck{S}\geq\cos\left(\frac{x^m}{2}\right)\left(x+\frac{x^m}{2}\right)-\sin\left(\frac{x^m}{2}\right),\\
\widecheck{S}\geq\frac{\sin\left({\underline{x}}\right)-\sin\left(\overline{x}\right)}{{\underline{x}-\overline{x}}}\left(x-{\underline{x}}\right)+\sin\left({\underline{x}}\right) \text{if~} \underline{x}\geq0,\\
\widecheck{S}\leq\frac{\sin\left({\underline{x}}\right)-\sin\left({\overline{x}}\right)}{{\underline{x}-\overline{x}}}\left(x-{\underline{x}}\right)+\sin\left({\underline{x}}\right) \text{if~} {\overline{x}}\leq0.
\end{cases}
\end{cases}\\
\label{eq:cosine envelope}
&\nonumber\left\langle\cos(x)\right\rangle^C=\\
&\quad\;\begin{cases}
\widecheck{C}:\begin{cases}
\widecheck{C}\leq 1-\frac{1 -\cos\left({x^m}\right)}{\left(x^m\right)^2}x^2,\\
\widecheck{C}\geq\frac{\cos\left(\underline{x}\right)-\cos\left({\overline{x}}\right)}{{\underline{x}-\overline{x}}}\left(x-{\underline{x}}\right)+\cos\left({\underline{x}}\right).  
\end{cases}
\end{cases}
\end{align}
\end{subequations}

\rn{where $x^m = \max\left(|\underline{x}|, |\overline{x}|\right)$. The auxiliary variables $\check{S}$ and $\check{C}$ represent the corresponding sets. For $-\frac{\pi}{2} < \underline{x} < \overline{x} < \frac{\pi}{2}$, the bounds on the sine and cosine functions are given as:}

\begin{subequations}
\begin{align}
& \underline{s} = \sin\left(\underline{x}\right) \leq \sin(x) \leq \overline{s} = \sin\left(\overline{x}\right), \\
\nonumber & \underline{c} = \min\left(\cos(\underline{x}),\cos(\overline{x})\right) \leq \cos(x) \\ & \quad \leq \overline{c} \!=\! \begin{cases} \max\left(\cos(\underline{x}),\cos(\overline{x})\right),\; \text{if~} \mathrm{sign}\left(\underline{x}\right) \!=\!  \mathrm{sign}\left(\overline{x}\right), \\ 1, \text{~otherwise}. \end{cases}\raisetag{1em}
\end{align}
\end{subequations} 

\rn{By replacing the squared, product, and trigonometric terms in Equations~\ref{eq:obj}--\ref{eq:qik1} with the auxiliary variables $w_{ii}$, $w_{lm}$, $c_{lm}$, and $s_{lm}$, the equations can be reformulated to incorporate convex envelopes. The resulting convexified equations are given below:}

\begin{subequations}
\label{eq:qc}
\begin{align}
&min\quad\sum_{i\in\mathcal{G}} c_{2,i}(P_{i}^{g})^2+c_{1,i}P_{i}^{g}+c_{0,i}\label{eq:obj3}\\
&\nonumber \text{subject to} \quad \left(\forall i\in\mathcal{B_A}, \forall   \left(l,m\right) \in\mathcal{L_A}\right)\\
\label{eq:qc_p}
& P_i^g-P_i^d = g_{sh,i}\, w_{ii}+\sum_{\substack{(l,m)\in \mathcal{L}\\ \text{s.t.} \hspace{3pt} l=i}} P_{lm}+\sum_{\substack{(l,m)\in \mathcal{L}\\ \text{s.t.} \hspace{3pt} m=i}} P_{ml}, \\
\label{eq:qc_q}
& Q_i^g-Q_i^d = -b_{sh,i}\, w_{ii}+\sum_{\substack{(l,m)\in \mathcal{L}\\ \text{s.t.} \hspace{3pt} l=i}} Q_{lm}+\sum_{\substack{(l,m)\in \mathcal{L}\\ \text{s.t.} \hspace{3pt} m=i}} Q_{ml},\\
\label{eq:qc_V}
&  (\underline{V}_i)^2\leq w_{ii} \leq (\overline{V}_i)^2,\\
\label{eq:qc_pik}
& P_{lm} = g_{lm} w_{ll} - g_{lm} c_{lm} - b_{lm} s_{lm}, \\
\label{eq:qc_qik}
& Q_{lm} = -\left(b_{lm}+b_{sh,lm}/2\right) w_{ii} + b_{lm} c_{lm}- g_{lm} s_{lm}, \\
\label{eq:qc_wii}
&  w_{ii} \in\left\langle V_i^2 \right\rangle^T, \\
\label{eq:qc_wik}
&  w_{lm} \in \left\langle V_l V_m \right\rangle^M, \\
\label{eq:qc_cik}
&  c_{lm} \in \left\langle w_{lm}\left\langle\cos\left(\theta_{lm} \right)\right\rangle^C\right\rangle^M, \\
\label{eq:qc_sik}
& s_{lm} \in \left\langle w_{lm} \left\langle\sin\left(\theta_{lm} \right)\right\rangle^S\right\rangle^M, \\
\label{eq:qc_jabr}
&  \text{Equations~}\eqref{eq:activelimit}\text{--}\eqref{eq:thetareflimit},\,\eqref{eq:cs_relationships_c},\,\eqref{eq:cs_relationships_s}.
\end{align}
\end{subequations}

\subsection{Customized QC-Branch and Bound Algorithm to solve OPF problem}

\rn{The accuracy of the QC relaxation method strongly depends on the quality of the variable bounds. Several studies have focused on Optimization-Based Bound Tightening techniques to improve these bounds in OPF problems~\cite{chen2015bound, coffrin2016strengthening}. These methods typically involve solving auxiliary optimization problems to iteratively tighten the upper and lower bounds of the variables. However, such approaches can become computationally intractable for large-scale systems due to their iterative nature.}
\rn{In this paper, we instead leverage the branch-and-bound (B\&B) algorithm, which systematically partitions the feasible space of the QC-relaxed OPF problem to search for improved solutions. The B\&B algorithm works by dividing the solution space into smaller subproblems (branching) and computing bounds on the optimal solution within each subproblem. By comparing these bounds, the algorithm can discard subproblems that cannot yield better solutions than the current best (pruning), thereby reducing the search space and improving computational efficiency. This method is particularly advantageous for non-convex problems like OPF, where finding the global optimum is especially challenging.}

\rn{To better understand the proposed method, we begin by reviewing the standard branch-and-bound algorithm. Consider a general optimization problem defined as \( \mathcal{O} = (\mathcal{D}, h) \), where \( \mathcal{D} \subseteq \mathbb{R}^n \) represents the set of feasible solutions, and \( h: \mathcal{D} \rightarrow \mathbb{R} \) is the objective function to be minimized. The aim is to identify an optimal solution \( \delta^\star \in \mathcal{D} \) such that, $\delta^\star \in \arg\min_{\delta \in \mathcal{D}} h(\delta)$. To find \( \delta^\star \), the Branch-and-Bound (B\&B) algorithm incrementally builds a search tree \( \mathcal{T} \), where each node (child) represents a subproblem defined over a subregion (branch) \( \mathcal{R} \subseteq \mathcal{D} \). The algorithm keeps track of the best solution found so far, called the \textit{incumbent}, and denotes it as \( \bar{\delta} \in \mathcal{D} \). At each step, it selects a subregion \( \mathcal{R} \in \Lambda \) to explore, where \( \Lambda \) is the current list of active subregions. If a new candidate solution \( \tilde{\delta} \in \mathcal{R} \) is found such that \( h(\tilde{\delta}) < h(\bar{\delta}) \), then the incumbent \( \bar{\delta} \) is updated to this better solution~\cite{}.}

\rn{At each step in the algorithm, there are two possibilities when exploring a region \( \mathcal{R} \): it can either be eliminated (pruned) or further divided for continued search.
If it is guaranteed that no feasible solution in \( \mathcal{R} \) can outperform the current best solution (the incumbent), that is,
$\forall \delta \in \mathcal{R}, \quad h(\delta) \geq h(\bar{\delta})$, then the region is \textit{pruned}, meaning it is excluded from further consideration. 
Otherwise, the region \( \mathcal{R} \) is divided into smaller subregions \( \mathcal{R}_1, \ldots, \mathcal{R}_m \), and these are added to the search tree \( \mathcal{T} \) for further exploration.
This process continues until every region has either been explored or pruned. As long as the pruning rules are correctly applied, the final incumbent \( \bar{\delta} \) is guaranteed to satisfy,
$\bar{\delta} \in \arg\min_{\delta \in \mathcal{D}} h(\delta)$.
}

\rn{The performance of the B\&B algorithm relies heavily on three strategic components including child selection strategy which determines the order in which subregions from \( \Lambda \) are explored, branching strategy which defines how each region \( \mathcal{R} \) is subdivided, and pruning strategy which determines when a region can be safely excluded. In many practical problems, such as in power systems, the domain \( \mathcal{D} \) is not listed explicitly. Instead, for any given candidate \( \delta \), we can quickly check whether it belongs to the domain, and we can generate subregions without needing to list every element in \( \mathcal{D} \).
The worst-case runtime of the B\&B algorithm depends on two factors: the branching factor \( \beta \), which is the maximum number of children each node can have, and the depth \( \zeta \) of the search tree \( \mathcal{T} \), which is how many times regions are divided. If solving a single subproblem takes at most \( \mu \) time, then the total time complexity is $O(\mu \beta^\zeta).$
However, pruning, one of the key features of the B\&B algorithm, greatly reduces this complexity in practice by eliminating regions that cannot contain a better solution.
}

\rn{In the context of power system optimization, particularly the non-convex AC-OPF problem, traditional nonlinear solvers often struggle to find globally optimal solutions. The B\&B framework addresses this challenge by partitioning the feasible space of teh problem and leveraging convex relaxations to compute lower bounds within each subregion, thereby facilitating a global solution to the non-convex OPF problem~\cite{gopalakrishnan2012global}.
B\&B algorithms, when combined with convex relaxation techniques, can solve nonconvex optimization problems to global optimality~\cite{sliwak2021semidefinite}. A broader overview of modern optimization techniques in power systems, including B\&B methods, is provided in~\cite{sun2017modern}. However, limited work has explored the integration of QC relaxations with B\&B for the OPF problem. This paper proposes a QC-based B\&B algorithm to address this gap.}








\rn{In the proposed QC-Relaxation B\&B approach, the objective function \( h(\theta) \) in the general B\&B formulation corresponds to the QC-OPF problem. The decision variable \( \eta \) includes the bus voltage magnitudes \( v = \{v_i\} \) and voltage angles \( \theta = \{\theta_i\} \). The resulting optimization problem is:
\[
\eta^\star = \{v^\star, \theta^\star\} \in \arg\min_{\eta \in \mathcal{D}} \text{QC-OPF}(\eta),
\]
subject to the feasible set \( \mathcal{D} \), defined by the OPF problem's operational and physical constraints.}

\rn{The proposed approach begins by solving a relaxed instance of the QC-OPF problem to initialize the root node of the B\&B tree.} 
\rn{At each level of the branching tree, we choose one variable to divide. This variable is either a bus voltage magnitude \( v_i \) within the range \( [v_i^{\min}, v_i^{\max}] \), or a voltage angle \( \theta_i \) within the range \( [\theta_i^{\min}, \theta_i^{\max}] \). In our method, we split a region \( \mathcal{R} \subseteq \mathcal{D} \) by dividing the selected variable's interval in half. That means we create two child regions by cutting the interval at its midpoint. For a variable \( z \in \{v_i, \delta_i\} \), the midpoint is calculated as:}
\[
z^{\text{mid}} = \frac{z^{\min} + z^{\max}}{2}.
\]

\rn{The left child keeps the same lower bound as the parent but sets its upper bound to the midpoint, that is,
\( z_{1} \in \left[ z^{\min},\ z^{\text{mid}} \right] \). The right child keeps the same upper bound as the parent but sets its lower bound to the midpoint:
\( z_{2} \in \left[ z^{\text{mid}},\ z^{\max} \right] \).
All other variables that are not chosen for branching keep their original bounds from the parent region.}


\rn{Each subregion, which we treat as a branch, defines a QC-OPF subproblem. This subproblem is solved to check whether a solution exists in that region (feasibility), and to find a local lower bound on the objective function. Specifically, for a subregion \( \mathcal{R} \subseteq \mathcal{D} \), the relaxed OPF problem is written as:}
\[
\eta^\mathcal{R} = \{v^\mathcal{R}, \theta^\mathcal{R}\} \in \arg\min_{\eta \in \mathcal{R}} \text{QC-OPF}(\eta),
\]
\rn{where \( \text{QC-OPF}(\eta) \) represents the QC relaxation of the original non-convex AC-OPF objective function over the region \( \mathcal{R} \).}

\rn{Relaxation methods, including the QC relaxation, can help determine if there is no feasible solution for the optimization problem within a specific region. We use this ability of the QC relaxation to eliminate subregions that have no feasible solutions. To do this, if the relaxed problem is found to be infeasible, that is, there is no \( \eta^\mathcal{R} \in \mathcal{R} \) that satisfies the convex constraints, then the region is discarded. 
Alternatively, let \( h_{\text{QC}}^\mathcal{R} \) be the value of the QC-OPF objective over region \( \mathcal{R} \), which is a lower bound for the OPF problem over region \( \mathcal{R} \), and let \( h(\bar{\eta}) \) be the current best (incumbent) objective value. If \( h_{\text{QC}}^\mathcal{R} \geq h(\bar{\eta}) \), then the region can also be safely discarded because it cannot provide a better solution than the incumbent.}
\rn{In addition, any region where the QC-OPF objective value \( h_{\text{QC}}^\mathcal{R} \) is greater than or equal to a valid upper bound from a feasible AC-OPF solution is also discarded. That is, if we have a valid upper bound \( h_{\text{AC}} \) from a feasible AC solution and:
\[
h_{\text{QC}}^\mathcal{R} \geq h_{\text{AC}},
\]
then the region \( \mathcal{R} \) is excluded from further consideration. On the other hand, if the region passes these checks, it is considered promising and kept for further splitting and evaluation in the next steps of the B\&B process~\cite{}.}

\begin{figure*}
    \centering
    \hspace{0.0cm}
\captionsetup{justification=centering}
\includegraphics[scale=0.48,trim= 12cm 0.20cm 10.7cm 0.17cm,clip]{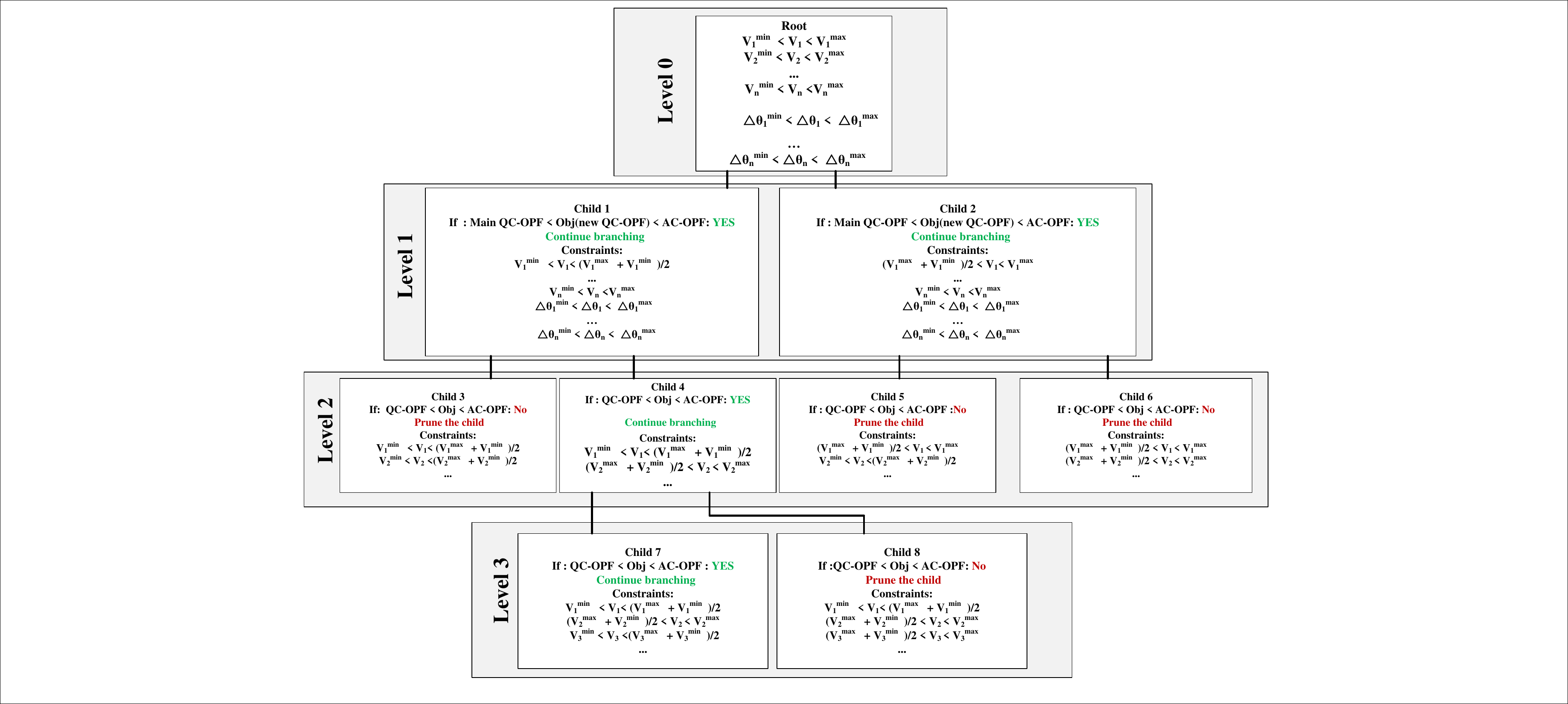}
	 \caption{\rn{Overview of the BB-QC-OPF framework. The flowchart illustrates how the Branch-and-Bound algorithm iteratively splits a single variable at each level, solving the QC-OPF problem within updated bounds for each child node. The figure highlights how variables' bounds are inherited from parent to child nodes throughout the branching process.}}
	\label{fig:general}
\end{figure*}

\rn{This process repeats until all subregions have been explored. By keeping track of a guaranteed global lower bound using the QC relaxation, and using a feasible AC-OPF solution as an upper bound, the algorithm is able to converge to the global optimal solution. In addition, the combination of refining intervals and systematically pruning unpromising regions makes the method both efficient and scalable. It effectively handles the non-convex nature of the AC-OPF problem while still ensuring high solution quality and global optimality. The steps of the proposed method are outlined in Algorithm~\ref{algorithm1}, and its overall process is illustrated in Figure~\ref{fig:general}.}


\begin{algorithm}
\caption{Branch-and-Bound for Optimal Power Flow with QC Relaxation}
\label{algorithm1}
\begin{algorithmic}[1]
\Require Network data $N$, QC-OPF bounds $\text{lower\_bound}$, $\text{upper\_bound}$, number of levels $n$, set of variables $\mathcal{V}$
\Ensure Best feasible solution $\text{best\_solution}$, QC lower bound $\text{QC\_lower\_bound}$, AC upper bound $\text{AC\_upper\_bound}$

\State \textbf{Initialization:}
\State Parse network data $N$ into initial QC-OPF relaxation (root node).
\State Solve QC-OPF to determine initial lower bound $\text{QC\_lower\_bound}$.
\State Solve AC-OPF to determine initial upper bound $\text{AC\_upper\_bound}$.
\State Initialize parent list $\mathcal{P}_0$ with root node.
\State $\text{best\_solution} \gets \emptyset$.

\For{$\ell = 1$ to $n$} \Comment{Iterative Branch-and-Bound Process}
    \State Initialize new parent list $\mathcal{P}_{\ell} \gets \emptyset$.
    \For{each $\text{parent} \in \mathcal{P}_{\ell-1}$}
        \State Select variable $v_{\ell} \in \mathcal{V}$ to split.
        \State Compute midpoint: $\text{midpoint} \gets \frac{\text{lower\_bound} + \text{upper\_bound}}{2}$.
        \State Generate two children:
        \State \hspace{1em} \textbf{Child 1:} Update $v_{\ell}^{\text{upper}} \gets \text{midpoint}$.
        \State \hspace{1em} \textbf{Child 2:} Update $v_{\ell}^{\text{lower}} \gets \text{midpoint}$.
        \For{each child}
            \State Solve QC-OPF for the child.
            \State Check feasibility.
            \State Compute objective value $f_{\text{child}}$.
            \If{$f_{\text{child}} < \text{QC\_lower\_bound}$ or $f_{\text{child}} > \text{AC\_upper\_bound}$ or infeasible}
                \State Prune the child.
            \Else
                \State Add the child to $\mathcal{P}_{\ell}$.
                \If{$f_{\text{child}} < f_{\text{best\_solution}}$}
                    \State Update $\text{best\_solution}$.
                \EndIf
            \EndIf
        \EndFor
    \EndFor
\EndFor

\State \Return $\text{best\_solution}$, $\text{QC\_lower\_bound}$, $\text{AC\_upper\_bound}$.

\end{algorithmic}
\end{algorithm}

\begin{figure}
    \centering
\captionsetup{justification=centering}
\includegraphics[scale=0.55]{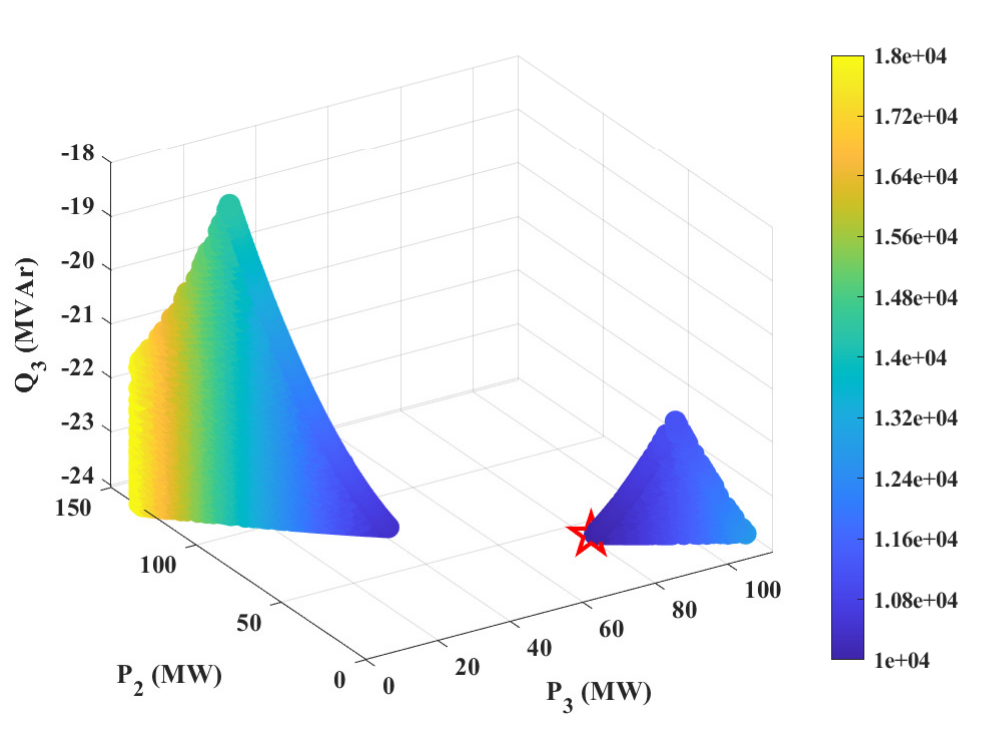}
	 \caption{\rn{Feasible space of the cyclic three bus system from~\cite{narimani2018empirical}, with $0.9 < V_1, V_2, V_3 < 1.1$ and $-2\pi < \Delta \theta_1, \Delta \theta_2, \text{and}~ \Delta \theta_3 < 2\pi$.}}
     
	\label{fig:sample0}
\end{figure}

\begin{figure*}[ht]
    \centering
    \begin{minipage}[t]{0.3\textwidth}
        \centering
        \includegraphics[width=\textwidth]{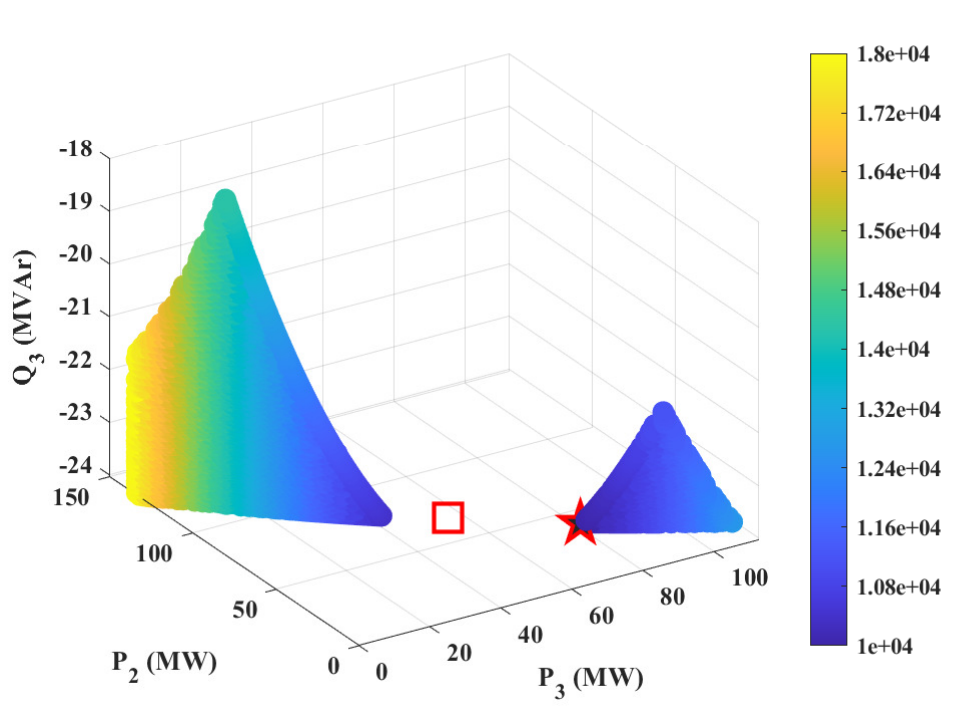}
        \caption{Feasible space of cyclic three bus system from~\cite{narimani2018empirical}, with $0.9<V_1<1$, $0.9<V_2<1.1$, $0.9<V_3<1.1$, and $-2\pi < \Delta \theta_1, \Delta \theta_2, \text{and} \Delta \theta_3 < 2\pi$.}
        \label{fig:sample1}
    \end{minipage}
    \hfill
    \begin{minipage}[t]{0.3\textwidth}
        \centering
        \includegraphics[width=\textwidth]{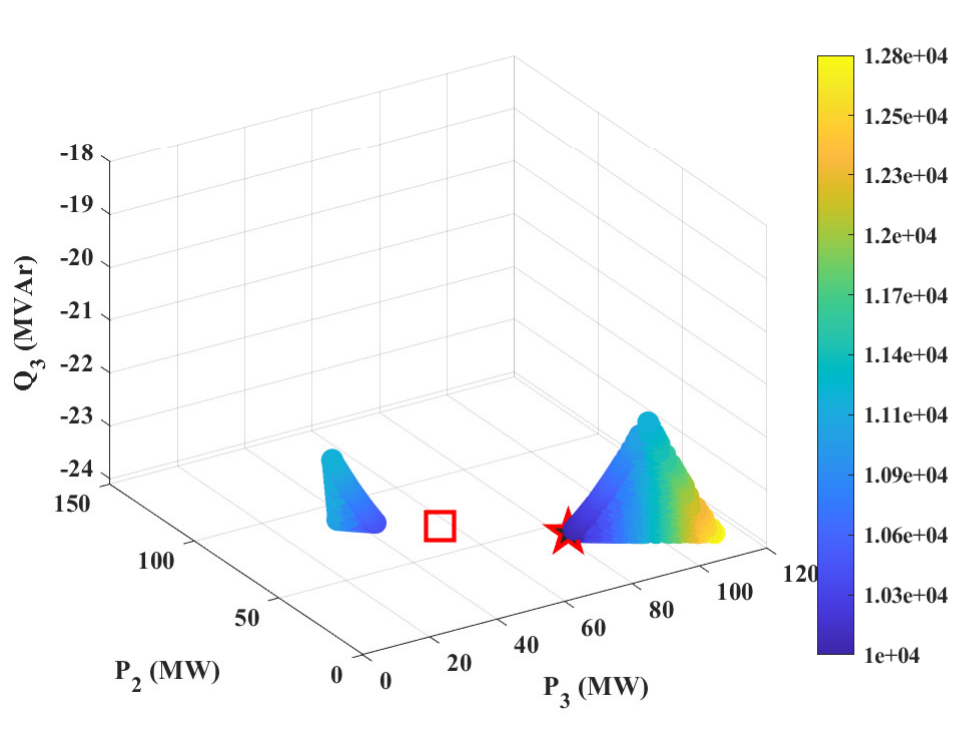}
        \caption{Feasible space of cyclic three bus system from~\cite{narimani2018empirical}, with $0.9<V_1<1$, $0.9<V_2<1$, $0.9<V_3<1.1$, and $-2\pi < \Delta \theta_1, \Delta \theta_2, \text{and} \Delta \theta_3 < 2\pi$.}
        \label{fig:sample3}
    \end{minipage}
    \hfill
    \begin{minipage}[t]{0.3\textwidth}
        \centering
        \includegraphics[width=\textwidth]{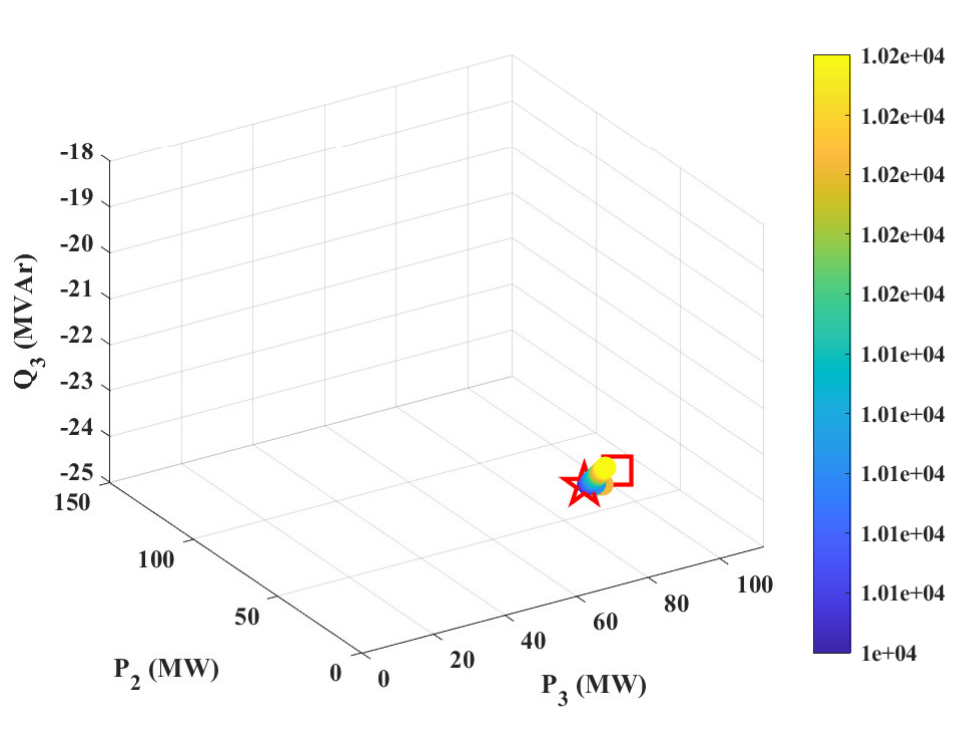}
        \caption{Feasible space of cyclic three bus system from~\cite{narimani2018empirical}, with $0.9<V_1<1$, $0.9<V_2<1$, $1<V_3<1.1$, and $-2\pi < \Delta \theta_1, \Delta \theta_2, \text{and}~ \Delta \theta_3 < 2\pi$.}
        \label{fig:sample7}
    \end{minipage}
\end{figure*}

\rn{To better understand the structure of the problem's feasible space, we analyzed the a three bus tests system in~\cite{}  and visualized how the proposed algorithm shapes its feasible region through repeated splitting and pruning. Specifically, we examined how the geometry of the OPF problem’s feasible space evolves after several steps of the proposed B\&B process.
Figure~\ref{fig:sample0} shows the overall feasible space of the OPF problem for this test system. In this example, six variables were considered: the voltage magnitudes at three buses and the angle differences across three branches, which were used as the splitting variables in the B\&B algorithm. Based on this setup, the resulting feasible spaces after multiple splitting and pruning steps are shown in Figs.~\ref{fig:sample1}, \ref{fig:sample3}, and \ref{fig:sample7}. From the visualizations, it becomes clear that after just a few iterations, the initially disconnected and complex feasible regions are reduced to a smaller, more connected region. This transformation makes it significantly easier for even a local solver to identify the global solution. This demonstrates how the proposed algorithm can effectively simplify the problem structure and guide the search toward the global optimum, even within a limited number of iterations.}

\section{Results and discussion}
\label{sec:results}

\rn{To test how well our proposed branch-and-bound method works for solving OPF problems, we applied it to several benchmark systems from the PGLib library \cite{babaeinejadsarookolaee2019power}. These systems are useful for testing because they are known to be challenging, there is often a big difference between the best-known solutions and the lower bounds from convex relaxation approaches.
We implemented our method using the Julia v1.11.5, JuMP v1.25.0~\cite{JuMP}, the PowerModels package~\cite{powermodels}, and the Gurobi 8.0 solver. The algorithm was designed to work smoothly within the PowerModels framework. It starts by extracting all necessary data, such as bus and branch information, voltage levels, and phase angles, from the case files provided by PowerModels.
We used Julia dictionaries to keep track of each step in the branch-and-bound process. These dictionaries store details about the subproblems (child nodes), their limits, and their objective values. As the algorithm runs, it updates these dictionaries to organize and store the results efficiently at every level of the method.
The results and analysis based on this setup are presented in the following sections.}

\subsection{Optimality Gaps and Solution Times}

\rn{In the first step, we applied the proposed method to several power system test cases. Table~\ref{tab:results_summary} provides a summary of these cases. They were chosen to show how the proposed approach can either improve the lower bound for the OPF problem or help find the global solution. In the table, the second column shows the objective value obtained from solving the AC OPF problem, which is used as the upper bound in the optimization. The following columns show how much the results from the QC-OPF and the proposed method differ from the AC-OPF results, expressed as a percentage gap. The last column shows the optimality gap, which is the difference between the local solution and the lower bound, normalized according to~\eqref{eq:optimality_gap}.}


\begin{equation}
\text{Optimality Gap} = \frac{\text{Local Solution - Lower Bound}}{\text{Local Solution}}
\label{eq:optimality_gap}
\end{equation}

\rn{The ``Local Solution'' refers to a result obtained by a local solver for the Optimal Power Flow (OPF) problem. In our approach, this acts as an upper bound. The ``Lower Bound'' in~\eqref{eq:optimality_gap} comes from either the QC relaxation or the proposed B\&B-QC method, both of which serve as lower bounds for the OPF problem. The optimality gap shows how close the lower bound from the B\&B-QC approach is to the AC-OPF solution. The remaining columns present the total time taken by the branch-and-bound algorithm for each case (in seconds), the number of branching levels, and the total number of child nodes created during the optimization.}

\rn{Upon comparison of the fourth and fifth columns in Table~\ref{tab:results_summary}, it is observed that the percentage gap between the AC and QC solutions is, in most instances, considerably larger than the gap between the AC and BB-QC solutions. This indicates that the QC-OPF solution can be refined by the branch-and-bound algorithm to more closely approximate the optimal AC-OPF value.
By splitting variables, the global solution of the OPF problem can be more readily obtained through the proposed algorithm. Nevertheless, a trade-off must be considered between reducing the optimality gap and the computational time required. The algorithm is affected by the curse of dimensionality, particularly in larger test cases. However, by applying the algorithm selectively to targeted variables, substantial improvements in the optimality gap can be achieved within a reasonable time frame.
Although the identification of the global solution is typically regarded as the primary goal of global optimization algorithms, the improvement of lower bounds is also of significant importance in various applications. As a key metric for evaluating solution quality, the objective value bounds derived from convex relaxations serve to indicate the proximity of a local solution to global optimality. Consequently, local algorithms and relaxations are frequently utilized in combination within spatial branch-and-bound frameworks for the solution of nonlinear programs (NLPs) and MINLPs~\cite{burer2012}.}

\begin{figure}
    \centering
\captionsetup{justification=centering}
\includegraphics[scale=0.5,trim= .5cm 0.60cm 1.0cm 0.17cm,clip]{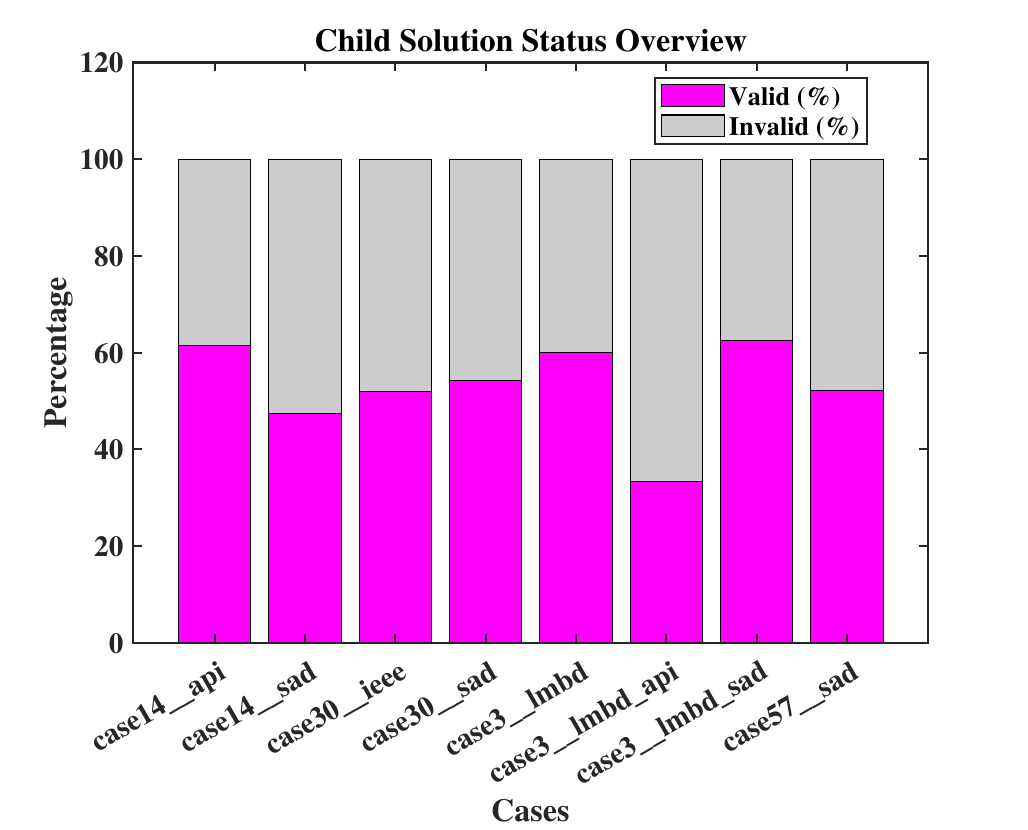}
	 \caption{\rn{Distribution of valid and invalid child solutions across different PGLib-OPF cases. The chart highlights the proportion of valid versus invalid children generated during the Branch-and-Bound process.}}
	\label{fig:column case}
\end{figure}

\subsection{\rn{Child Solution Status Overview}}

\rn{To assess the pruning efficiency and decision-making effectiveness of the proposed branching algorithm, we analyze the distribution of valid and invalid child nodes across multiple test cases. Fig.~\ref{fig:column case} presents the percentage of valid and invalid child nodes generated during the branching process for various test cases. Each bar corresponds to a specific case and is divided into two segments: the magenta segment represents the percentage of valid child nodes, those satisfying all problem constraints, while the gray segment indicates the percentage of invalid child nodes, which are either infeasible or violate problem bounds and are therefore pruned.
This statistical overview highlights a significant proportion of child nodes that are invalid across different cases. The figure underscores the effectiveness of the proposed QC relaxation technique, which enables early identification and pruning of non-promising child nodes. By comparing each child node’s objective value with a high-quality lower bound derived from the QC relaxation, the algorithm is able to discard unproductive branches early in the search process. This reduces the computational burden and accelerates convergence without sacrificing solution quality.
Furthermore, the results demonstrate that an effective branching strategy, when combined with strong relaxations, can substantially shrink the search space by eliminating many unviable subproblems. Consequently, the proposed approach enhances both the efficiency and robustness of the solution process in solving non-convex optimization problems such as the AC-OPF.}

\begin{table*}
\centering
\caption{Summary of optimization results for different cases.}
\begin{tabular}{|l|c|c|c|c|c|c|}
\hline
\textbf{Case Name} & \textbf{AC Obj} & \thead{\textbf{QC}\\\textbf{Gap (\%)}} & \thead{\textbf{BB-QC}\\\textbf{Gap (\%)}} & \thead{\textbf{Execution}\\\textbf{Time (Sec)}} & \thead{\textbf{$N_{\mathrm{Lev}}$}} & \thead{\textbf{$N_{\mathrm{Child}}$}} \\ \hline
case3\_lmbd              & 5812.64         & 0.98                 & 0.1                     & 0.93                            & 6                         & 22                            \\ \hline
case3\_lmbd\_\_api       & 11242.12        & 4.79                 & 2.99                    & 0.81                            & 6                         & 12                            \\ \hline
case3\_lmbd\_\_sad       & 5959.31         & 1.4                  & 0.09                    & 1.04                            & 6                         & 20                            \\ \hline
case14\_ieee             & 2178.08         & 0.11                 & 0.0004                  & 81.13                           & 34                        & 852                           \\ \hline
case14\_ieee\_\_api      & 5999.36         & 5.4                  & 0.02                    & 593.04                          & 34                        & 5974                          \\ \hline
case14\_ieee\_\_sad      & 2776.78         & 23.69                & 3.89                    & 5.51                            & 34                        & 38                            \\ \hline
case24\_ieee\_rts\_\_sad & 76917.96        & 2.81                 & 0                       & 270.02                          & 62                        & 1216                          \\ \hline
case30\_ieee             & 8208.51         & 22.96                & 0.01                    & 1205.53                         & 71                        & 4742                          \\ \hline
case30\_ieee\_\_api      & 18036.58        & 5.73                 & 0.02                    & 17587.27                        & 71                        & 17062                         \\ \hline
case30\_ieee\_\_sad      & 8208.51         & 6.01                 & 0.07                    & 1139.56                         & 30                        & 5050                          \\ \hline
case57\_ieee\_\_sad      & 38663.28        & 0.32                 & 0                       & 11277.64                        & 57                        & 29090                         \\ \hline
\end{tabular}
\label{tab:results_summary}
\end{table*}

\subsection{Branching Process Analysis}

\begin{figure}
    \centering
\captionsetup{justification=centering}
\includegraphics[scale=0.59]{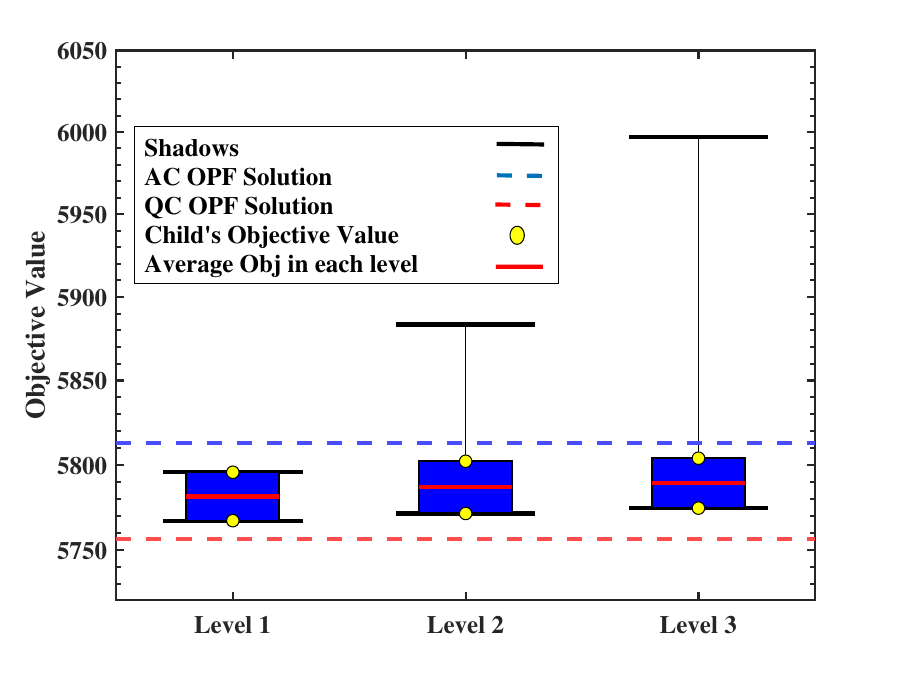}
	 \caption{\rn{Candle chart illustrating the objective value distribution across three branching levels for ``Case3\_lmbd''. The chart includes the floor (best child objective), roof (worst valid child objective), upper shadow (gap to QC-OPF), and average value (mean of child objectives) at each level, shown relative to the AC-OPF and QC-OPF benchmarks.}}
	\label{fig:candle3bus}
\end{figure}

\begin{figure}
    \centering
\captionsetup{justification=centering}
\includegraphics[scale=0.53]{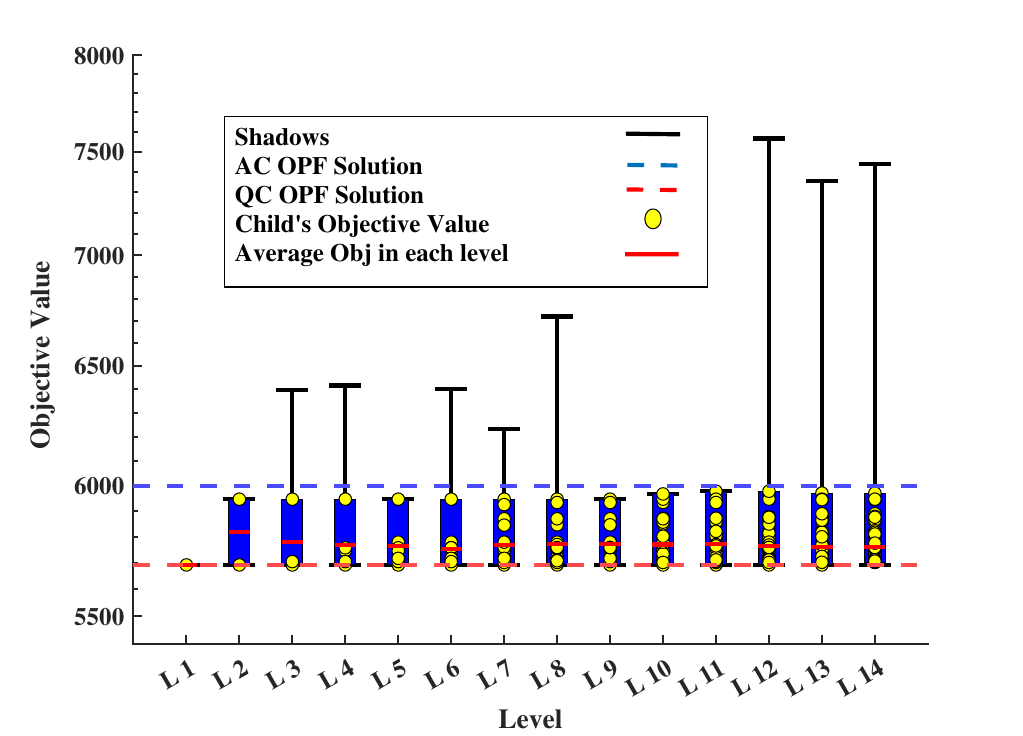}
	 \caption{\rn{Candle chart illustrating the objective value distribution across three branching levels for ``Case14\_ieee''. The chart includes the floor (best child objective), roof (worst valid child objective), upper shadow (gap to QC-OPF), and average value (mean of child objectives) at each level, shown relative to the AC-OPF and QC-OPF benchmarks.}}
	\label{fig:candle14bus}
\end{figure}

\begin{figure}
    \centering
\captionsetup{justification=centering}
\includegraphics[scale=0.5,trim= 1.7cm 0.08cm 1.0cm 0.17cm,clip]{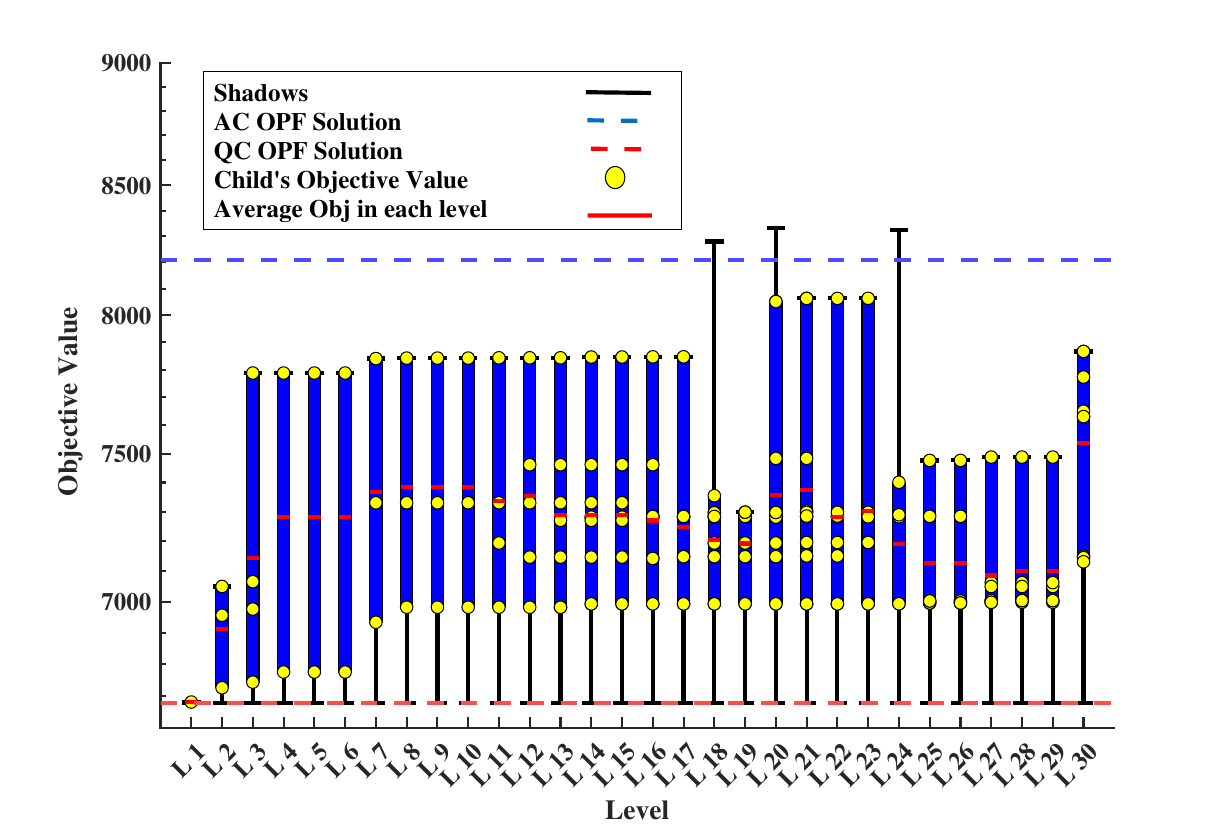}
	 \caption{\rn{Candle chart illustrating the objective value distribution across three branching levels for ``Case30\_ieee''. The chart includes the floor (best child objective), roof (worst valid child objective), upper shadow (gap to QC-OPF), and average value (mean of child objectives) at each level, shown relative to the AC-OPF and QC-OPF benchmarks.}}
	\label{fig:candle30bus}
\end{figure}

\begin{figure*}
    \centering
\captionsetup{justification=centering}
\includegraphics[scale=0.64,trim= 1.5cm 0.20cm 1.0cm 0.17cm,clip]{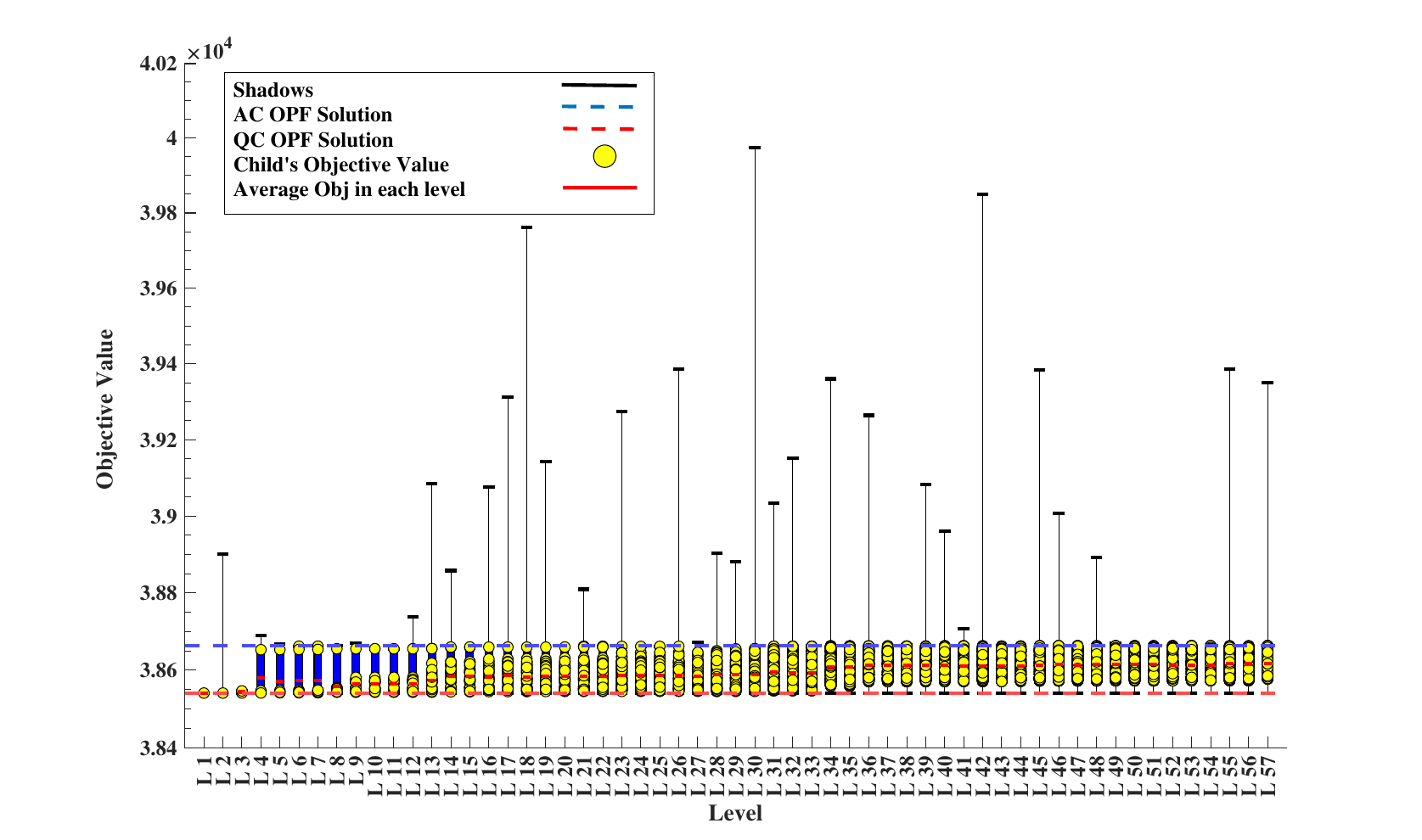}
	 \caption{\rn{Candle chart illustrating the objective value distribution across three branching levels for ``Case57\_ieee''. The chart includes the floor (best child objective), roof (worst valid child objective), upper shadow (gap to QC-OPF), and average value (mean of child objectives) at each level, shown relative to the AC-OPF and QC-OPF benchmarks.}}
	\label{fig:candle57bus}
\end{figure*}

\rn{To better understand the dynamics of the branching process and evaluate the performance of the proposed algorithm, we present a set of candle charts in Figs. \ref{fig:candle3bus} through \ref{fig:candle57bus}. These visualizations illustrate the objective value distributions of child nodes at each branching level for different test cases. In this analysis, branching decisions are based solely on bus voltage variables, and the candle charts offer a hierarchical view of how the algorithm explores and prunes the solution space.}
\rn{Each candle chart summarizes key statistical features of the objective values at each level of the branching tree. The upper shadow represents the maximum objective value among all children at that level, regardless of feasibility. The lower shadow indicates the minimum objective value. The roof and floor represent the maximum and minimum objective values, respectively, among valid child nodes, i.e., those satisfying all problem constraints. A red horizontal line within each candle marks the mean of the valid objective values, providing an indicator of solution quality at that level. Yellow dots represent the individual objective values of all child nodes, enabling a visual assessment of the distribution and density of explored solutions. In addition, two reference lines are included: the blue dashed line shows the AC-OPF objective value, while the red dashed line represents the QC-OPF relaxation bound.
This representation allows us to assess how effectively the algorithm prunes invalid or suboptimal child nodes and whether the branching decisions lead to the generation of high-quality solutions. Narrow candles, where roofs and floors are close, suggest a concentrated set of valid solutions and an effective local search around promising areas. Wide candles indicate more variance in solution quality and potentially weaker branching performance at those levels. The length of the shadows, particularly the difference between the upper shadow and the roof, reveals the extent of infeasible or non-competitive children that were generated and subsequently pruned. Furthermore, the closeness of the average line and the roof to the AC solution is a strong indicator of algorithmic performance: smaller gaps imply that the branching strategy is successfully guiding the search toward globally optimal regions, while larger gaps may signal inefficiencies or opportunities for improving branching rules.}

\rn{Fig.~\ref{fig:candle3bus} depicts the case for ``Case3-lmbd'', covering three levels of branching. The candles in this case are narrow, with both the average line and the roof remaining close to the AC-OPF solution. This indicates that most valid children at each level have objective values near the global optimum and that the algorithm effectively identifies and concentrates on promising regions early in the search. The relatively short shadows reflect the limited presence of poor-quality or infeasible children, confirming that the pruning strategy, guided by the QC lower bound, is working efficiently in the early stages of this smaller case.
Fig.~\ref{fig:candle14bus} extends the analysis to a deeper tree for ``Case14-ieee'', spanning 14 levels. In this scenario, we observe more variation across levels. While the earlier levels exhibit moderately narrow candles and close alignment with the AC objective, some mid-levels introduce broader candles and slightly lower average values. This behavior suggests an initial phase of focused exploration followed by a broader search in later levels. Nonetheless, the proximity of the average lines to the AC bound remains satisfactory in most levels, reflecting that the algorithm maintains overall solution quality and continues to prune suboptimal branches effectively.}

\rn{Fig.~\ref{fig:candle30bus} analyzes the more complex ``Case30-IEEE'' scenario over 30 levels. This case shows increased variability both in terms of candle width and shadow height. In early levels (e.g., L1–L10), the chart reveals wider candles and greater shadow lengths, indicating a mix of valid and invalid children and a broader search of the solution space. However, as the levels progress, the average lines stabilize and gradually approach the AC objective line, suggesting that the algorithm is able to recover and guide the search toward optimal regions. This behavior demonstrates the algorithm’s adaptive nature, adjusting from an exploratory phase into more refined pruning and selection as the search deepens.
Fig.~\ref{fig:candle57bus} presents the branching process for ``Case57-IEEE'', the largest test case. Despite the size, the candle chart reveals consistent performance across all 57 levels. Most candles remain narrow with average values closely tracking the AC-OPF solution. This consistency highlights the scalability and robustness of the proposed branching strategy. While long upper shadows appear in certain levels, indicating the generation of poor-quality children, these are efficiently pruned, and they do not adversely impact the distribution of valid solutions. The clustering of yellow points within each level further reinforces that the algorithm is effectively concentrating its search efforts around high-quality regions even at greater depth.
These candle chart analyses demonstrate that the proposed branching method, supported by QC relaxation bounds and voltage-based decision rules, effectively balances exploration and pruning. It consistently identifies and focuses on promising subregions of the solution space, particularly evident in smaller cases like ``Case3-lmbd'' and large-scale scenarios like ``Case57-IEEE''. The results confirm the method’s capacity to deliver high-quality solutions with improved computational efficiency, offering strong potential for solving large non-convex optimization problems such as AC-OPF.}

\subsection{\rn{Impact of Child Nodes on Optimality Gap Improvement}}

\rn{This section investigates how child nodes generated during the proposed B\&B process contribute to narrowing the optimality gap between the QC relaxation and the solution form local solver, i.e., AC OPF solution. Specifically, we aim to understand how the distribution of normalized objective values across children reflects the effectiveness of branching in tightening lower bounds and advancing toward the global optimum. By analyzing the shape and progression of these distributions, we gain insight into the intrinsic difficulty of each test case, the quality of the initial QC relaxation, and the efficacy of variable branching decisions.}

\begin{figure*}
    \centering
\captionsetup{justification=centering}
\includegraphics[scale=0.605]{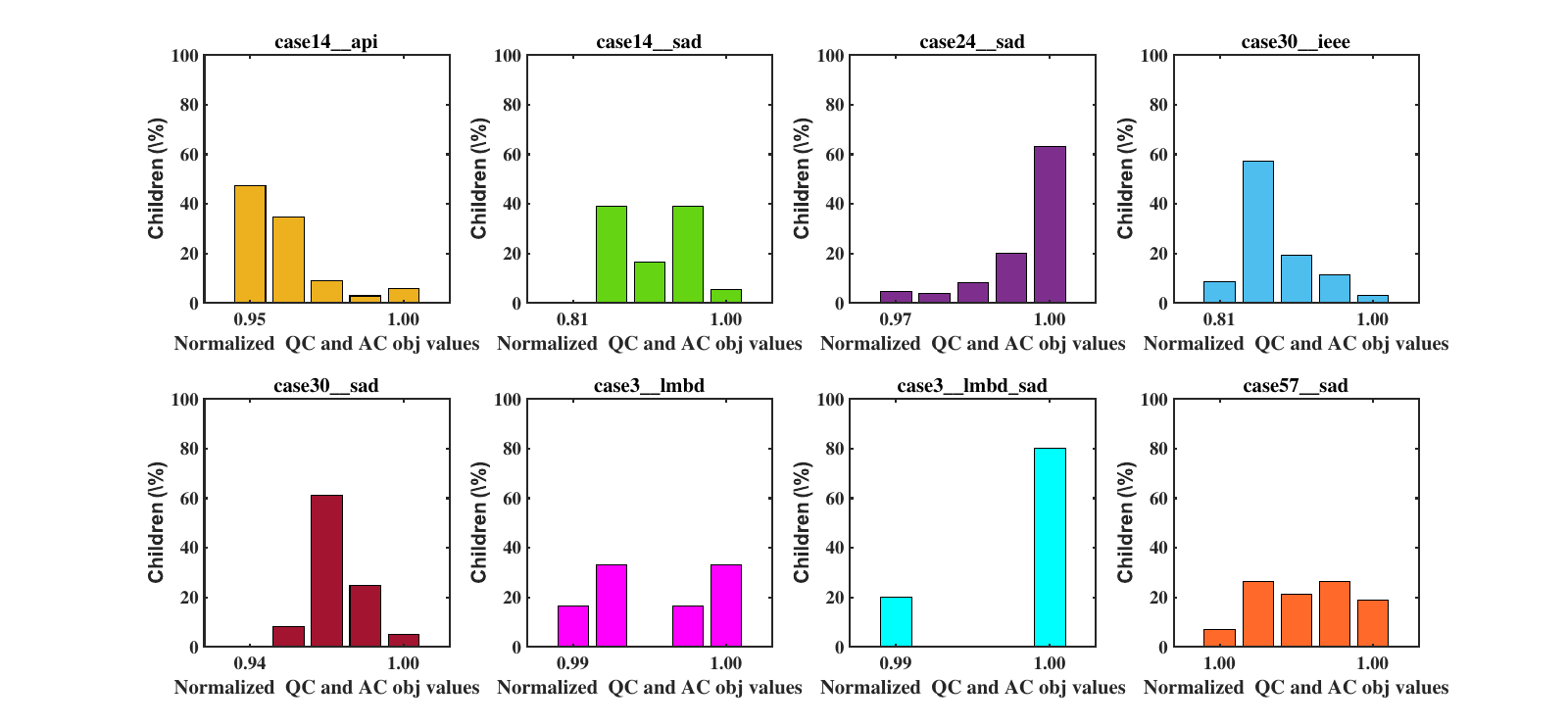}
	 \caption{\rn{Histograms showing the distribution of normalized objective values for child solutions across different test cases. Each chart illustrates how the proposed method improves solution quality relative to the QC and AC objective values, highlighting the progression toward optimality within each case.}}
	\label{fig:histogram}
\end{figure*}

\rn{Figure~\ref{fig:histogram} presents histograms for a set of benchmark cases, each showing the distribution of child-node objective values normalized by the AC-OPF optimum:
\[
  \text{Normalized Cost} = \frac{OPF_\mathrm{child}}{OPF_\mathrm{AC}}.
\]
The horizontal axis spans from the QC relaxation lower bound \(\bigl(OPF_\mathrm{QC}/OPF_\mathrm{AC}\bigr)\) on the left to the AC solution \(\bigl(OPF_\mathrm{AC}/OPF_\mathrm{AC}\bigr)\) at \(1.00\) on the right. This interval is subdivided into five equal-width bins. Each bar represents the percentage of generated children falling within a given normalized objective range.}

\rn{Initially, all child nodes inherit the QC bound and thus populate the leftmost bin. As branching progresses, tighter relaxations are achieved through variable splits, causing child objective values to increase and spread toward the AC solution. This rightward movement of mass across the histogram indicates successful improvement of the lower bound and provides a visual trace of how the B\&B algorithm navigates the feasible region. For instance, in the ```case14\_api'' scenario, the QC relaxation underestimates the AC cost by approximately 5.4\%, resulting in a histogram that begins at a normalized value of roughly \(0.946\). The histogram reveals that over 80\% of child nodes remain within the first two bins near the QC bound. Although some progress is made, with at least one child falling in the rightmost bin, the small spread emphasizes the difficulty of achieving substantial improvement early in the branching process. This suggests that the initial branching yields limited improvement in the objective value, but as branching progresses through multiple iterations, the objective function gradually approaches the AC solution more closely. By contrast, ``case14\_sad'' exhibits a more evenly distributed histogram: 40\% of children fall in the second bin, 17\% in the third, and the remainder spread across the rest. Despite a wider initial gap (23.69\%), the histogram reflects a relaxation that improves more effectively with branching, likely due to better alignment between the QC model and the problem’s nonlinear structure.}

\rn{In ``case24\_sad'', a dramatic shift occurs, with only 9\% of children in the first two bins and 64\% in the right-most bin. This concentration near the AC optimum suggests that the feasible space is readily exploitable and that early branching steps rapidly tighten the bounds. In contrast, the two ``case30'' variants show a different pattern, many of the early iterations fall into bins near the QC bound, with only a few late-branching solutions reaching the right-most bins. This suggests that while the initial branching steps effectively tighten the relaxation by cutting into the feasible space, this tightening does not substantially advance the solution toward the AC optimum until the final iterations.}

\rn{The ``case3\_lmbd'' problem exhibits a histogram where early iterations show limited impact on closing the optimality gap, followed by a sudden shift after a specific branching iteration, where the solution quickly approaches the AC solution and remains in that vicinity through subsequent iterations. This behavior indicates that meaningful improvements in optimality emerge only after sufficient branching depth is reached. Finally, ``case57\_sad'' produces a nearly uniform distribution, revealing a moderate relaxation that necessitates multiple stages of refinement across various regions of the feasible set. This pattern indicates that most iterations contribute relatively evenly to improving the solution, with each branching step gradually tightening the relaxation and incrementally reducing the optimality gap.}

\rn{These histograms offer a visual and quantitative summary of each case's ``branching'' process. A histogram skewed to the left indicates a loose initial relaxation, with most children still far from the AC optimum, thus requiring extensive exploration. A right-heavy distribution, on the other hand, reflects rapid convergence enabled by effective relaxations and well-chosen branching decisions.
Furthermore, since the B\&B framework employed here splits one variable at a time, shifts in the histogram can also identify which variables (e.g., voltage magnitudes or angle differences at specific buses) contribute most to improving the optimality gap. This analysis motivates enhancements such as prioritizing impactful variables in the branching rule or incorporating valid inequalities related to those variables to guide the B\&B more efficiently through challenging regions.
Beyond qualitative assessment, these findings also have practical implications. Cases exhibiting rapid rightward progression are more computationally tractable and may benefit from early termination criteria or aggressive node pruning. In contrast, cases with left-heavy or flat distributions may require deeper lookahead strategies, adaptive branching rules, or hybrid relaxations to make meaningful progress. By combining histogram-based diagnostics with dynamic strategy selection, future B\&B implementations could more effectively allocate computational resources and improve scalability for large-scale power systems.}

\subsection{\rn{Comparative Analysis of Optimality Gaps Using Candlestick Visualization}}

\begin{figure}
    \centering
\captionsetup{justification=centering}
\includegraphics[scale=0.42,trim= 1.7cm 0.10cm 3.0cm 1.27cm,clip]{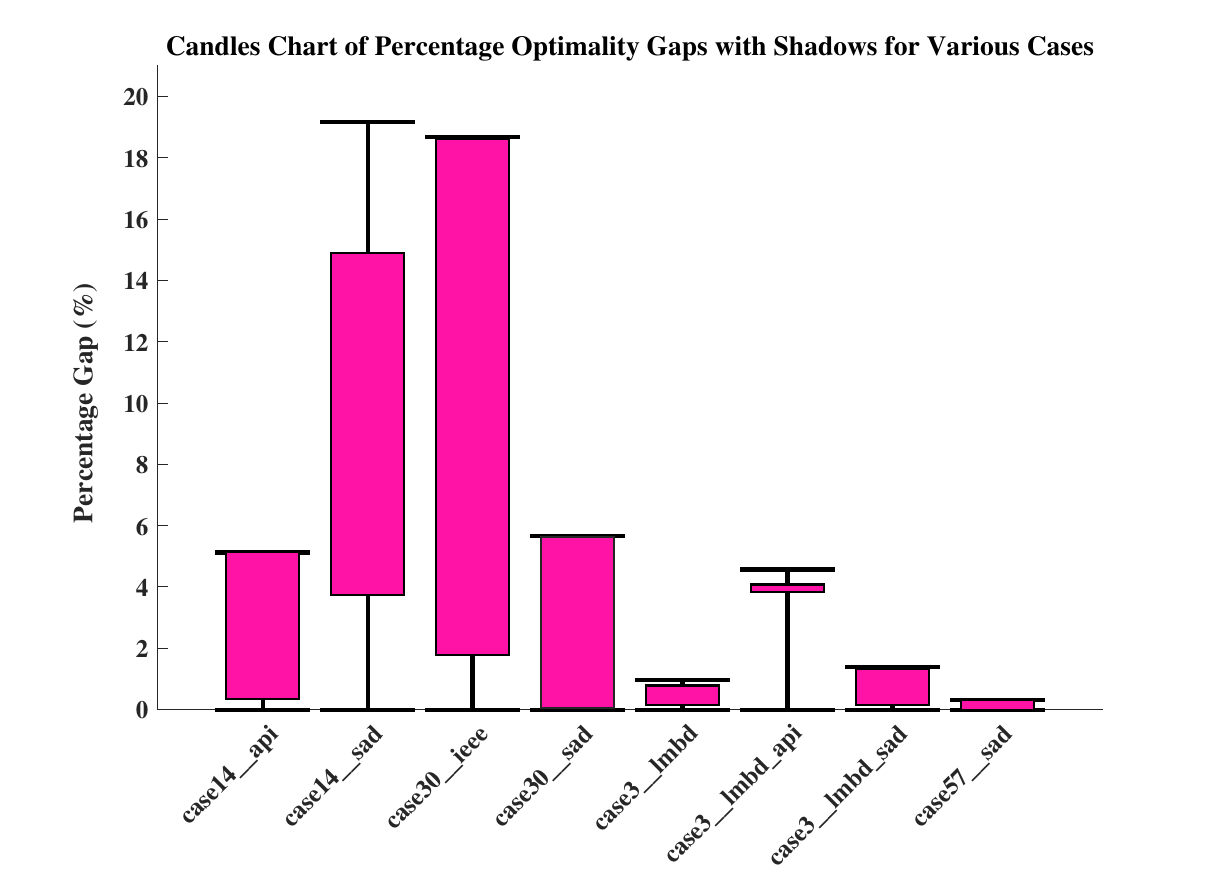}
	 \caption{\rn{Optimality gap distributions for selected PGLib-OPF cases using the proposed QC-assisted Branch-and-Bound method. Candlesticks show solution spread and relaxation strength relative to the AC-OPF benchmark.}}
	\label{fig:allcandles}
\end{figure}

\rn{To further evaluate the effectiveness of the proposed branching strategy, we conduct a comparative analysis of optimality gaps across multiple test cases. This assessment uses a modified candlestick chart to visualize the distribution and progression of solution quality relative to the AC-OPF benchmark. Figure~\ref{fig:allcandles} provides a visual summary of the optimality gaps across several power system test cases using a modified candlestick chart format. This chart compares the quality of solutions obtained through the proposed branching strategy with respect to the AC-OPF benchmark. Each candlestick captures the range of percentage gaps between the AC-OPF objective value and three critical metrics: the best-performing child solution, the worst-performing valid child solution, and the initial QC relaxation. Specifically, the bottom of the candle, referred to as the floor, represents the percentage gap between the AC-OPF value and the best child solution found. The top of the candle, or roof, denotes the gap to the worst valid child solution. The upper shadow reflects the gap between the initial QC relaxation and the AC-OPF solution, while the lower shadow is zero, as it represents the difference between the AC solution and itself.}

\rn{An analysis of this chart reveals notable differences in optimization performance across the test cases. For instance, ``case14\_sad'' and ``case30\_ieee'' display tall candlesticks with prominent upper shadows, indicating that the initial QC relaxations were relatively weak, leaving substantial room for improvement. The wide candle bodies in these cases further suggest high variability among child solutions, pointing to a complex and potentially fragmented feasible region. These cases require deeper and more targeted branching to navigate the solution space effectively and to approach the AC-OPF benchmark more closely. Although the presence of some near-optimal child solutions is promising, the overall variability highlights the need for refinement in either the relaxation technique or the branching strategy. In contrast, ``case3\_lmbd'' shows a much shorter candlestick, with a relatively small gap between the best child solution and the AC solution. This behavior implies that while early branching steps have limited effect, the optimization process significantly improves once a certain depth is reached, after which the solution stabilizes near the AC benchmark. This sudden and sustained improvement reflects a relaxation that becomes significantly more effective after the feasible space is sufficiently refined. Likewise, ``case3\_lmbd\_api'' and ``case3\_lmbd\_sad'' present very narrow candles with minimal spread, suggesting that the branching strategy performs consistently well and converges quickly due to a relatively well-structured feasible region.}

\rn{The most balanced performance is observed in ``case57\_sad'', which features a uniformly short candlestick. This indicates that improvements in the optimality gap are achieved gradually and consistently across iterations. There is no sudden leap in performance; rather, each branching step incrementally tightens the relaxation and improves solution quality. This behavior suggests that the branching strategy is particularly effective for this instance, offering a well-calibrated balance between feasible space exploration and convergence. The candlestick chart in Figure~\ref{fig:allcandles} serves as a diagnostic tool for evaluating the performance of the branching strategy across diverse cases. By examining the shape and spread of each candle, one can assess the quality of the relaxation, the variability of child solutions, and the convergence behavior of the method. These insights are essential for identifying cases that require stronger relaxations, improved branching heuristics, or more refined solution techniques, ultimately contributing to more robust and efficient algorithms for nonconvex power system optimization.}

\section{Conclusion}
\label{sec:Conclusion}

\rn{In this paper, we proposed a QC-assisted B\&B algorithm for solving the OPF problem. The approach integrates QC relaxation techniques, which are effective in providing lower bounds for nonconvex problems, with the conventional B\&B framework to systematically divide the OPF feasible space. By combining the global search capabilities of B\&B with the bounding strength of QC relaxations, the method improves both solution accuracy and computational efficiency. The algorithm operates by iteratively branching on selected problem variables, including voltage magnitudes and phase angle differences, and choosing branches that yield the most promising improvements in the objective function. This strategy allows the method to effectively reduce the search space while guiding the optimization process toward globally optimal or near-optimal solutions. Additionally, the algorithm provides insights into the impact of branching decisions on individual variables, helping to identify how different bounds and splits influence convergence and solution quality. We applied the proposed method to various cases from the standard PGLib-OPF test suite to evaluate its performance. The results demonstrate the algorithm’s capability to consistently find global or near-global solutions with improved computational performance compared to traditional methods. The proposed approach successfully combines the abilities of convex relaxation and global optimization techniques, and lays the groundwork for future research in efficient, scalable, and accurate OPF solution methods.}

\bibliographystyle{IEEEtran}
\IEEEtriggeratref{50}
\bibliography{ref}

\end{document}

\end{document}

\section{Conclusion}
\label{sec:conclusion}

This paper has delved deeply into the realm of AC False Data Injection (FDI) attacks within power systems, aiming to unravel their complexities and devise effective countermeasures. To address these challenges, we have meticulously delineated the constraints across different segments of the attack region, guided by the diverse physical constraints in each segment. By adopting the principle of "When in Rome, do as the Romans do," we have tailored and examined FDI attack designs on the IEEE 39 and 118 test bus systems in such a way that BDD could not be able to detect them while the state variables are deviated from their real values. However, although this paper has tried to demonstrate a unified approach to designing AC FDI attack, but the attack vector for implementing this proposed attack is not optimum. Over-shooting of residuals in the SE process indicates this subject. In fact, this proposed attack can be implemented with an optimum scenario about the number of the measurements for attacking, or values of the attack vector to reduce the over-shootings in the residuals of estimations, or finding an optimum zone or multi-zone to run a successful AC attack. These points serve as potential avenues for future research in AC attack designs, which may commence with the unified approach proposed in this paper and subsequently undergo optimization across diverse scenarios and methodologies.

\end{document}